\documentclass{article}
\usepackage{graphicx}
\usepackage{amsmath}
\usepackage{amsfonts}
\usepackage{amssymb}
\newtheorem{theorem}{Theorem}

\newtheorem{condition}{Condition}

\newtheorem{corollary}{Corollary}

\newtheorem{definition}{Definition}
\newtheorem{example}{Example}

\newtheorem{lemma}{Lemma}
\newtheorem{notation}{Notation}

\newtheorem{proposition}{Proposition}
\newtheorem{remark}{Remark}

\newenvironment{proof}[1][Proof]{\textbf{#1.} }{\ \rule{0.5em}{0.5em}}

\begin{document}

\title{On the determination of 2-step solvable Lie algebras from its weight graph}
\author{Jos\'{e} Mar\'{\i}a Ancochea\thanks{corresponding author : e-mail : Jose\_Ancochea@mat.ucm.es} \and Otto Rutwig
Campoamor\thanks{%
research supported by the D.G.Y.C.I.T project PB98-0758} \\
%EndAName
Departamento de Geometr\'{\i}a y Topolog\'{\i}a\\
Fac. CC. Matem\'{a}ticas Univ. Complutense\\
28040 Madrid ( Spain )}
\date{}
\maketitle

\begin{abstract}
By using the concept of weight graph associated to nilpotent Lie algebras $\frak{g}$, we find necessary and sufficient conditions for a semidirect product $\frak{g}\oplus T_{i}$, where $T_{i}<T$ is a subalgebra of a maximal torus of derivations $T$ of $\frak{g}$ which induces a decomposition of $\frak{g}$ into one dimensional weight spaces, to be 2-step solvable. In particular we show that the semidirect product of such a Lie algebra with its torus of derivations cannot be itself 2-step solvable.\newline AMS Subj. Clas. 05C99, 17B30\newline \textit{keywords : Lie algebra, nilpotent, weight, graph, rigid}    
\end{abstract}

\section{Introduction}
Graph theory has been proven to be useful for the study of Lie algebras, as follows for example from the classical theory of semisimple Lie algebras and its Dynkin diagrams. Though it is well known that these diagrams only exist for this particular class of Lie algebras, the main idea of associating graphs can be applied to study weight systems of nilpotent Lie algebras. However, one important assumption such as the one dimensionality of weight spaces gets lost, which forces us to distinguish two cases : those algebras admitting a decomposition into one dimensional weight spaces ( for the action of an exterior torus of derivations ), and those not. For the first class ( simple ) graphs can be used, while for the latter we have to consider di- and pseudographs.\newline
In this work we associate to a nilpotent Lie algebra without multiple weights a so called weight graph, which describes the structure of the weight system associated to the Lie algebra. In particular, this graph is of interest when studying the 2-step solvable subalgebras of the semidirect product $\frak{g}\oplus T$, when $T$ is a maximal torus of derivations over $\frak{g}$. \newline These results about weight graphs can also be used to establish an interesting result about the rigidity of such algebras, property that otherwise would be considerably complicated to formulate, as for dimensions $n\geq 8$ there do not exist a classification of nilpotent Lie algebras.
\newline For elementary properties and definitions about Lie algebras, we refer to reference [4], while for questions about graph theory we refer to [5]. Here we shall only consider finite undirected graphs without loops or multiple edges.

\subsection{Preliminaries and notations}

Let $\frak{g}$ be a finite dimensional complex nonsplit nilpotent Lie algebra. Let
$Der\left(  \frak{g}\right)  $ be its Lie algebra of derivations. A torus $T$
over $\frak{g}$ is an abelian subalgebra of $Der\left(  \frak{g}\right)  $
consisting of semi-simple endomorphisms. Clearly the torus $T$ induces a
natural representation [4] on the Lie algebra $\frak{g}$, such that this
decomposes as
\[
\frak{g}=\bigsqcup_{\alpha\in T^{\ast}}\frak{g}_{\alpha}%
\]
where $T^{\ast}=Hom_{\mathbb{C}}\left(  T,\mathbb{C}\right)  $ and
$\frak{g}_{\alpha}=\left\{  X\in\frak{g\;}|\;\left[  t,X\right]
=\alpha\left(  t\right)  X\;\;\forall t\in T\right\}  $ is the weight space
corresponding to the weight $\alpha$. If the torus is maximal for the
inclusion relation, as tori are conjugated [6], its common dimension is a
numerical invariant of $\frak{g}$ called the rank and denoted by $r\left(
\frak{g}\right)$. Following Favre [3], we call
\[
R\frak{g}\left(  T\right)  =\left\{  \alpha\in T^{\ast}\;|\;\frak{g}_{\alpha
}\neq0\right\}
\]
the set of weights for the representation of $T$ over $\frak{g}$ and
\[
P\frak{g}\left(  T\right)  =\left\{  \left(  \alpha,d\alpha\right)
\;|\;\alpha\in R\frak{g}\left(  T\right)  ,\;d\alpha=\dim\frak{g}_{\alpha
}\right\}
\]

\begin{definition}
Let $\frak{g}$ be a nilpotent Lie algebra and $T$ a maximal torus of
derivations. Then $R\frak{g}\left(  T\right)  $ is called a weight system for
$\frak{g}$.
\end{definition}

\begin{remark}
Two sets $P\frak{g}\left( T\right)$ and $P\frak{g}\left( T'\right)$ are called equivalent 
if they induce an equivalent representation over the nilpotent Lie algebra $\frak{g}$.
The equivalence class of a weight system constitutes an invariant of the algebra
[3].
\end{remark}

Though the weights can be chosen quite arbitrarily, the
relations they satisfy are preserved.\newline Given the $n$-dimensional Lie algebra $\frak{g}$, we have $R\frak{g}\left(T\right)  =\left\{  \alpha_{1},..,\alpha_{r}\right\}  $, where $r\leq n$. In order to avoid the use of pseudo- or digraphs, here we will only consider the case $n=r$, which allows us to identify both $P\frak{g}\left(T\right)$ and $R\frak{g}\left(T\right)$.
For certain results we will have to add an additional condition :

\begin{condition}
There exists a minimal presentation $\left<X,R\right>$ of the algebra $\frak{g}$ such that the linear system $S$ associated to this presentation [1] is equivalent to the linear system defined by the weights of $R\frak{g}\left(T\right)$.
\end{condition}

\begin{remark}
The condition implies that for the chosen representative of the isomorphism class of $\frak{g}$ a sum of weights which is also a weight is induced by some structural constants.\newline 
It is not difficult to prove that all nilpotent Lie algebras of maximal nilindex satisfy the preceding condition. Other wide classes satisfying it can be found in [1].
\end{remark}

Suppose that $R\frak{g}\left(  T\right)  =\left\{  \alpha_{1},..,\alpha
_{n}\right\}  $ is the weight system of $\frak{g}$. To this weight system we
can associate the following graph : let $V\left(  G\right)  =\left\{
v_{1},..,v_{n}\right\}  $ be the points, where $v_{i}$ corresponds to the weight $\alpha_{i}$ for all $i$.
We say that $\alpha_{i}$ joins $\alpha_{j}$ if $\alpha_{i}+\alpha_{j}\in
R\frak{g}\left(  T\right)  $. The corresponding graph is denoted by $G\left(
R\frak{g}\left(  T\right)  \right)  $.

\begin{lemma}
Let $R\frak{g}\left(  T\right)  $ be a weight system for the nilpotent Lie
algebra $\frak{g}$. Then $G\left(  R\frak{g}\left(  T\right)  \right)  $
contains at least an isolated point.
\end{lemma}

\begin{proof}
As the algebra is finite dimensional, it follows from the weight space
decomposition that there exists at least one weight $\gamma\in R\frak{g}%
\left(  T\right)  $ which is ''extreme'' in the following sense :%
\[
\alpha+\gamma\notin R\frak{g}\left(  T\right)  ,\;\forall \alpha\in R\frak{g}%
\left(  T\right)
\]
Then for any vector $X\in\frak{g}_{\gamma}$ we have $ad_{\frak{g}}X\equiv0$,
where $ad_{\frak{g}}$ denotes the adjoint operator in $\frak{g}$. Thus the
vector is central, and the corresponding vertex of $G\left(  R\frak{g}\left(
T\right)  \right)  $ is isolated.
\end{proof}

\begin{corollary}
The graph $G\left(  R\frak{g}\left(  T\right)  \right)  $ contains $\left(
\dim Z\left(  \frak{g}\right)  \right)  $ isolated points, where $Z\left(
\frak{g}\right)  $ is the center of $\frak{g}$.
\end{corollary}

\subsection{The weight graph $G=\overline{G\left(  R\frak{g}%
\left(  T\right)  \right)  }$}

We denote by $\overline{G\left(  R\frak{g}\left(  T\right)
\right)  \text{ }}$ the complementary graph to $G\left(  R\frak{g}\left(
T\right)  \right)  $. It follows immediately from the lemma that
$\overline{G\left(  R\frak{g}\left(  T\right)  \right)  }$ is a connected graph.

\begin{definition}
Let $R\frak{g}\left(  T\right)  $ be a weight system of $\frak{g}$. Then the
graph $\overline{G\left(  R\frak{g}\left(  T\right)  \right)  }$ is called the
weight graph of $\frak{g}$.
\end{definition}

\begin{lemma}
For the complementary $G$ of a $\left(  p,q\right)$-weight graph $\overline{G}$ we have
\[
q\leq\sum_{j=1}^{\left[  \frac{p}{2}\right]  }\left(  p-2j\right)
\]
where $\left[  \frac{p}{2}\right]  $ denotes the integer part of $\frac{p}{2}$.
\end{lemma}

\begin{proof}
As we have imposed $\dim\frak{g}_{\alpha}=1$ for any weight $\alpha\in
R\frak{g}\left(  T\right)  $, we can reorder the weights of $\frak{g}$ in such
manner that a relations $\alpha_{i}+\alpha_{j}=\alpha_{k}$ corresponds to a
sum $i+j=k$, where $k\leq p$. Thus the maximal number of sums of weights
equals the number of possibilities $i+j=k$ with $k\leq p$ and $1\leq i<j$. It
is easily seen that this number is precisely $\sum_{j}^{\left[  \frac{\pi}%
{2}\right]  }\left(  p-2j\right)  $.
\end{proof}

\begin{remark}
Observe that this bound has been obtained indepently of graph theory. It is
only based in the principal property of weights, namely, that certain sums of
them give another weight.
\end{remark}

\begin{proposition}
Let $G$ be a $\left(  p,q\right)  $-weight graph of a nilpotent Lie algebra
$\frak{g}$. Then
\[
q\geq\left(
\begin{array}
[c]{c}%
p\\
2
\end{array}
\right)  -\sum_{j=1}^{\left[  \frac{p}{2}\right]  }\left(  p-2j\right)
\]
\end{proposition}

\begin{lemma}
For $p\geq 4$ the following inequality holds :
\begin{equation*}
\left( 
\begin{array}{c}
p-1 \\ 
2
\end{array}
\right) -\sum_{i=1}^{\left[ \frac{p}{2}\right] }\left( p-2j\right) >0
\end{equation*}
\end{lemma}

\begin{proof}
For $p=4$ the assertion is obvious. Now suppose it holds for $p>4$. Then 
\begin{eqnarray*}
\left( 
\begin{array}{c}
p \\ 
2
\end{array}
\right) -\sum_{i=1}^{\left[ \frac{p+1}{2}\right] }\left( p+1-2j\right) 
&=&\left( 
\begin{array}{c}
p \\ 
2
\end{array}
\right) -\sum_{i=1}^{\left[ \frac{p}{2}\right] }\left( p-2j\right) -\left[ 
\frac{p}{2}\right]  \\
&=&\left( 
\begin{array}{c}
p-1 \\ 
2
\end{array}
\right) -\sum_{i=1}^{\left[ \frac{p+1}{2}\right] }\left( p-2j\right) +\left(
p-1\right) -\left[ \frac{p}{2}\right] >0
\end{eqnarray*}
\end{proof}

Recall that a graph $G$ is called bipartite if there exists a partition of
its set of points, $V\left( G\right) =V_{1}\bigcup V_{2}$, such that each
line in $G$ joins a point of $V_{1}$ with a point of $V_{2}$. 

\begin{proposition}
A $\left( p,q\right) $-weight graph $\overline{G\left( R\frak{g}\left(
T\right) \right) }$ is not bipartite.
\end{proposition}

\begin{proof}
Suppose that $\overline{G\left( R\frak{g}\left( T\right) \right) }$ is
bipartite. As it has a point of degree $\left( p-1\right) $, say $v_{0}$,
the partition of the point set $V$ must be $V_{1}=\{v_{0}\},\;V_{2}=V-V_{1}$.
This implies that there are no lines in the weight graph joining points of $%
V_{2}$, which implies that the subgraph whose point set is $V_{2}$ is
totally disconnected. Taking the complementary, it follows that $G\left( R%
\frak{g}\left( T\right) \right) $ has complete subgraph of degree $\left(
p-1\right) $, isomorphic to $K_{p-1}$. Now this graph has $\left( 
\begin{array}{c}
p-1 \\ 
2
\end{array}
\right) $ lines, but $G\left( R\frak{g}\left( T\right) \right) $ has at most 
$\sum_{i=1}^{\left[ \frac{p-1}{2}\right] }\left( p-2j\right) $ lines. The
contradiction follows from the preceding lemma.
\end{proof}

\begin{corollary}
Any weight graph has at least one odd cycle.
\end{corollary}

In fact we can improve this by showing that, up to the lowest dimensional case, all weight graphs have triangles in common :

\begin{proposition}
For $p\geq 4$ any weight graph contains a triangle.
\end{proposition}

\begin{proof}
The only case to be considered is when the weight graph $\overline{G\left( R\frak{g}\left( T\right) \right) }$ is associated to a nilpotent Lie algebra whose center is one dimensional. For any higher dimensional centers, let $v_{1},v_{2}$ be two vertices corresponding to central weight spaces, and $v_{3}$ another arbitrary point. In the complementary $G\left( R\frak{g}\left( T\right) \right)$ these points are totally disconnected, thus they induce a triangle in the weight graph. So we can assume that the graph $G\left( R\frak{g}\left( T\right) \right)$ has a unique isolated point. If the weight graph does not contain a triangle, then $G\left( R\frak{g}\left( T\right) \right)$ contains a complete subgraph isomorphic to $K_{p-1}$, which is not possible in view of prop. 2.
\end{proof}

\begin{corollary}
A weight graph $\overline{G\left( R\frak{g}\left( T\right) \right) }$ is a tree if and only if $\frak{g}$ is isomorphic to the three dimensional Heisenberg Lie algebra.
\end{corollary}

\begin{corollary}
If $r\left(\frak{g}\right)\geq 3$, then any maximal complete subgraph of $\overline{G\left( R\frak{g}\left(T\right) \right) }$ has degree $r\leq p-2$.
\end{corollary}

\section{The fundamental subgraphs $\overline{G\left( \protect\beta
_{i}\right) }$ }

Suppose that $\frak{g}$ has rank $k$ and that $T$ is a maximal torus. Then
we can write $T=\sum T_{i}$, where $T_{i}$ $\left( 1\leq i\leq k\right) $ is
the one dimensional subtorus defined by its weights
\begin{eqnarray*}
\lambda _{i}^{j} &=&\delta _{ij}\beta _{i},\;1\leq j\leq k \\
\lambda _{i}^{k+j} &=&a_{i}^{j}\beta _{i},\;1\leq j\leq \dim \,\frak{g}-k
\end{eqnarray*}
and $a_{i}^{j}\in \mathbb{C}$. Then we can write the weight set as 
\begin{equation*}
R\frak{g}\left( T\right) =\left\{ \beta _{1},..,\beta _{k},\gamma
_{k+1},..,\gamma _{n-k}\right\} 
\end{equation*}
where 
\begin{equation*}
\gamma _{k+j}=\sum_{i=1}^{k}a_{i}^{j}\beta _{i},\;1\leq j\leq \dim \frak{g}-k
\end{equation*}
Let $n=dim \frak{g}$. Then we can find a basis $\left\{ X_{1},..,X_{n}\right\} $ of $\frak{g}$
such that 
\begin{eqnarray}
\frak{g} &=&\frak{n}^{\beta _{j}}+..+\frak{n}^{\beta _{k}}+\frak{n}^{\gamma
_{1}}+..\frak{n}^{\gamma _{n-k}} \\
&=&\mathbb{C}X_{1}+..+\mathbb{C}X_{n}  \notag
\end{eqnarray}
We denote the semidirect product $\frak{g}\oplus T_{i}$ of $\frak{g}$ and
the torus $T_{i}$ $\left( 1\leq i\leq k\right) $ by $T\left( \beta
_{i}\right) $, and for $1\leq s\leq k$ the semidirect product 
$\frak{g}\oplus T_{i_{1}}\oplus ..\oplus T_{i_{s}}$ by $T\left( \beta
_{i_{1}},..,\beta _{i_{s}}\right) $. 

\begin{definition}
Let $1\leq i\leq k$. For a weight $\gamma \in R\frak{g}\left( T\right) $ the 
$\beta _{i}$-length of $\gamma $ is defined as 
\begin{equation*}
l_{\beta _{i}}\left( \gamma \right) =a_{i}
\end{equation*}
where $\gamma =\sum_{i=1}^{k}a_{i}\beta _{i}$.
\end{definition}

\begin{notation}
For $1\leq i\leq k$ we denote by $E\left( \beta _{i}\right) $ the set 
\begin{equation*}
E\left( \beta _{i}\right) =\left\{ \gamma \in R\frak{g}\left( T\right)
\;|\;l_{\beta _{i}}\left( \gamma \right) \geq 1\right\} 
\end{equation*}
\end{notation}

\begin{definition}
For $1\leq i\leq k$ let $\overline{G\left( \beta _{i}\right) }$ be the
subgraph of $\overline{G\left( R\frak{g}\left( T\right) \right) }$ whose
point set is $E\left( \beta _{i}\right) $. The subgraph $\overline{G\left(
\beta _{i}\right) }$ is called the fundamental subgraph associated to the (
fundamental ) weight $\beta _{i}$.
\end{definition}

\begin{remark}
It can easily happen that two nonequivalent weight systems have the same
weight graph. However, the structure of its fundamental subgraphs shows
their different structure.
\end{remark}

\begin{theorem}
A weight system $R\frak{g}\left( T\right) $ $\ $is completely determined by
its weight graph $\overline{G\left( R\frak{g}\left( T\right) \right) }$ $\ $%
and its fundamental subgraphs $\overline{G\left( \beta _{i}\right) }$. 
\end{theorem}

\begin{proof}
Let $R\frak{g}\left( T\right) =\left\{ \beta _{1},..,\beta _{k},\gamma
_{k+1},..,\gamma _{n-k}\right\} $ and $R\frak{g}^{\prime }\left( T^{\prime
}\right) =\left\{ \beta _{1}^{\prime },..,\beta _{k}^{\prime },\gamma
_{k+1}^{\prime },..,\gamma _{n-k}^{\prime }\right\} $ be two nonequivalent
weight systems associated to nilpotent Lie algbras of rank $k$, where 
\begin{equation*}
\gamma _{k+j}=\sum_{i=1}^{k}a_{i}^{j}\beta _{i};\;\gamma _{k+j}^{\prime
}=\sum_{i=1}^{k}b_{i}^{j}\beta _{i}^{\prime }\; 1\leq j\leq n-k
\end{equation*}
Suppose that both their associated weight graphs $\overline{G_{1}}$ and $%
\overline{G_{2}}$ and its fundamental subgraphs $\overline{G_{1}\left( \beta
_{i}\right) }$ and $\overline{G_{2}\left( \beta _{i}^{\prime }\right) }$ \ $%
\left( 1\leq i\leq k\right) $ are isomorphic. Then we can find a permutation 
$\sigma \in S_{k}$ such that 
\begin{equation*}
\sigma \left( G_{1}\left( \beta _{i}\right) \right) =G_{2}\left( \beta
_{i}^{\prime }\right) ;\;1\leq i\leq k
\end{equation*}
for the complementary of the fundamental subgraphs.\newline
Define the map $f:R\frak{g}\left( T\right) \rightarrow R\frak{g}^{\prime
}\left( T^{\prime }\right) $%
\begin{equation*}
\beta _{i}\mapsto \beta _{i}^{\prime }
\end{equation*}
Then, whenever the weight $\beta _{i}$ is joined with $\beta _{j}$, its
image $\beta _{i}^{\prime }$ is joined with $\beta _{j}^{\prime }$, as it
corresponds to a line of the subgraph $\overline{\overline{G_{1}\left( \beta
_{i}\right) }}$. It follows at once that  $R\frak{g}\left( T\right) $ and $R%
\frak{g}^{\prime }\left( T^{\prime }\right) $ are the same weight system,
which contradicts the hipothesis.
\end{proof}

\begin{remark}
Of course, the underlying nilpotent Lie algebras need not to be isomorphic.
Although there is only a finite number of equivalence classes for weight
systems [3], for dimensions $n\geq 7$ there is an infinity of isomorphism
classes for nilpotent Lie algebras.
\end{remark}

Let $\frak{r}$ be a solvable non-nilpotent Lie algebra. Then $\frak{r}$ is called $s$-step solvable if $D^{\left( s\right) }\frak{r}=0$ and $D^{\left( s-1\right) }\frak{r}\neq 0$,
where $D^{\left( k+1\right) }\frak{r}=\left[ D^{\left( k\right) }\frak{r}%
,D^{\left( k\right) }\frak{r}\right] $ are the terms of the derived sequence.

\begin{lemma}
Let $\overline{G\left( \beta _{i}\right) }$ be a fundamental subgraph such
that there exist two weights $\gamma _{1},\gamma _{2}\in E\left( \beta
_{i}\right) $ such that $\gamma _{1}+\gamma _{2}\in R\frak{g}\left( T\right) 
$ we have $\gamma _{1}+\gamma _{2}\notin E\left( \beta _{i}\right) $. Then 
the semidirect product $\frak{g}\oplus T_{i}$ is at least three step
solvable.
\end{lemma}

The proof is trivial. Observe that the sum of the two weights corresponds to
a zero weight for the action of the torus $T_{i}$ defined by its weights $E\left(\beta_{i}\right)$, thus these points are not connected in the weight
graph.\newline
From this property we easily obtain a characterization of 2-step solvable semidirect products $\frak{g}\oplus T_{i}$ : 

\begin{theorem}
Let $\frak{g}$ be a $n$-dimensional nilpotent Lie algebra of rank $k$. Then
the semidirect product $T\left( \beta _{i}\right) $ is 2-step solvable if
and only if the fundamental subgraph $\overline{G\left( \beta _{i}\right) }$
is complete.
\end{theorem}

\begin{proof}
Let $\left\{ X_{1},..,X_{n}\right\} $ be a basis as in $\left( 1\right) $. As the
semidirect product $T\left( \beta _{i}\right) $ is 2-step solvable, whenever
there are vectors $X_{i},X_{j}$ such that $\left[ X_{i},X_{j}\right] \neq 0$%
, we have $X_{i}\in \frak{n}^{\gamma }$ with $\gamma \notin E\left( \beta
_{i}\right) $ or $X_{j}\in \frak{n}^{\gamma }$ with $\gamma \notin E\left(
\beta _{i}\right) $. Thus for any vector $X\in \frak{n}^{\gamma }$ with $%
\gamma \in E\left( \beta _{i}\right) $ we have 
\begin{equation*}
\left[ X,\sum_{\alpha \in E\left( \beta _{i}\right) }\frak{n}^{\alpha }%
\right] =0
\end{equation*}
This implies, in view of the construction of the graph $G\left( R\frak{g}%
\left( T\right) \right) $,  that the complementary to the fundamental
subgraph $\overline{G\left( \beta _{i}\right) }$ is totally disconnected, so
that $\overline{G\left( \beta _{i}\right) }$ is isomorphic to the graph $%
K_{p}$, where $p=card\,E\left( \beta _{i}\right) $. \newline
The converse follows at once.
\end{proof}

This results gives in fact the key to decide whether a semidirect product $%
T\left( \beta _{i_{1}},..,\beta _{i_{s}}\right) $ is 2-step solvable or not.
Let us denote by $\overline{G\left( \beta _{i_{1}},..,\beta _{i_{s}}\right) }
$ the subgraph of $\overline{G\left( R\frak{g}\left( T\right) \right) }$
whose points belong to $\bigcup_{j=1}^{s}E\left( \beta _{i_{j}}\right) $.
Clearly this is the subgraph obtained by removing the points which are not
in the union. 

\begin{corollary}
Let $2\leq s\leq k$. Then the semidirect product $T\left( \beta
_{i_{1}},..,\beta _{i_{s}}\right) $ is 2-step solvable if and only if the
subgraph $\overline{G\left( \beta _{i_{1}},..,\beta _{i_{s}}\right) }$ is
complete.
\end{corollary}

\begin{example}
Let $\frak{g}$ be the Lie algebra with presentation
\[
\left< X_{1},..,X_{6}| [X_{1},X_{i}]=X_{i+1}, (i=2,3), [X_{2},X_{5}]=-[X_{3},X_{4}]=X_{6}\right>
\]
It is of rank $3$ and its weight system is 
\[
\left\{\beta_{1},\beta_{2},\beta_{3},\beta_{1}+\beta_{2},\beta_{1}+\beta_{3},\beta_{1}+\beta_{2}+\beta_{3}\right\}
\]
The graph $G\left( R\frak{g}\left( T\right) \right)$ is isomorphic to a nonramified tree with an isolated point. Clearly all of its fundamental subgraphs $G\left( \beta _{i}\right)$ are totally disconnected, thus the corresponding semidirect products are 2-step solvable. It can be easily seen that there does not exist a two dimensional subtorus with this property.
\end{example}

We call the join of two graphs $G_{1}$ and $G_{2}$
with disjoint point \ and line sets, $V_{1}$ and $V_{2}$, respectively $L_{1}
$ and $L_{2}$, to the graph $G_{1}+G_{2}$ consisting of the union $G_{1}\cup
G_{2}$ ( in the obvious sense ) and all lines joining $V_{1}$ and $V_{2}$. 

\begin{proposition}
Let $\overline{G\left( \beta _{i_{1}},..,\beta _{i_{s}}\right) }$ be a
complete subgraph of $\overline{G\left( R\frak{g}\left( T\right) \right) }$.
If $E\left( \beta _{i_{r}}\right) \cap E\left( \beta _{i_{t}}\right)
=\emptyset $ for $r\neq t$, then   $\overline{G\left( \beta
_{i_{1}},..,\beta _{i_{s}}\right) }$ is isomorphic to the join  $\overline{%
G\left( \beta _{i_{1}}\right) }+..+$ $\overline{G\left( \beta
_{i_{s}}\right) }$. 
\end{proposition}

\begin{proof}
Let $\overline{G\left( \beta _{i_{1}},..,\beta _{i_{s}}\right) }$ be
complete and fix a weight, $\beta _{i_{1}}$ for example. Then the removal of
the set $\bigcup_{2\leq j\leq s}E\left( \beta _{i_{j}}\right) $ from the
weight graph gives the fundamental subgraph $\overline{G\left( \beta
_{i_{1}}\right) }$, as the intersections of the sets $E\left( \beta
_{i}\right) $ are empty. Thus 
\begin{equation*}
\overline{G\left( \beta _{i_{1}},..,\beta _{i_{s}}\right) }=\overline{%
G\left( \beta _{i_{1}}\right) }+\overline{G\left( \beta _{i_{2}},..,\beta
_{i_{s}}\right) }
\end{equation*}
By recurrence it follows that 
\begin{equation*}
\overline{G\left( \beta _{i_{1}},..,\beta _{i_{s}}\right) }=\overline{%
G\left( \beta _{i_{1}}\right) }+...+\overline{G\left( \beta _{i_{s}}\right) }
\end{equation*}
\end{proof}

From this result we deduce an important property of the algebras $\frak{g}\oplus T$ :

\begin{proposition}
Let $\frak{g}$ be a nonabelian nilpotent Lie algebra and $T$ a maximal torus
of derivations. Then the semidirect product $\frak{g}\oplus T$ is at least
three step solvable.
\end{proposition}

\begin{proof}
If the semidirect product were 2-step solvable, then its weight graph would
be complete. Then its complementary graph $G\left( R\frak{g}\left( T\right)
\right) $ $\ $is totally disconnected, which implies that there are no
nontrivial brackets in $\frak{g}.$
\end{proof}

\section{Applications to rigid Lie algebras}

In this section we apply the information obtained about the fundamental subgraphs of a weight graph $\overline{G\left( R\frak{g}\left( T\right) \right) }$ to prove that there do not exist  2-step solvable rigid Lie algebras. For this purpose we recall the elementary facts about rigidity.
\bigskip Let $\mathcal{L}^{n}$ be the algebraic variety of complex Lie algebra laws on
$\mathbb{C}^{n}$. Consider the natural action of the algebraic group
$GL\left(  n,\mathbb{C}\right)  $ on $\mathcal{L}^{n}$ given by
\begin{align*}
GL\left(  n,\mathbb{C}\right)  \times L^{n}  &  \rightarrow L^{n}\\
(f,\mu)  &  \rightarrow f*\mu
\end{align*}
with $f*\mu(X,Y)=f^{-1}(\mu(f(X),f(Y))$, for all $X,Y\in\mathbb{C}^{n}.$ We
note by $\mathcal{O(\mu)}$ the orbit of $\mu$.

\begin{definition}
The Lie algebra law $\mu$ (or the complex Lie algebra $\frak{g}$ of law $\mu)$
is called rigid if $\mathcal{O(\mu)}$ is a Zariski open subset in $\mathcal{L}^{n}.$
\end{definition}

Each open orbit of this natural action of $GL\left(  n,\mathbb{C}\right)  $ on
$\mathcal{L}^{n}$ gives, considering its Zarisky closure, an irreducible component of
$\mathcal{L}^{n}$ .\ Therefore, only a finite number of those orbits exists; or,
equivalently, there is only a finite number of isomorphism classes of rigid Lie
algebras with open orbit.
\newline Now it is known that any rigid Lie algebra is decomposable, i. e., it can be written as 
\[
\frak{g}=\frak{s}\oplus\frak{t}\oplus\frak{n}%
\]
where $\frak{s}$ is a Levi subalgebra, $\frak{n}$ the nilradical and
$\frak{t}$ an abelian subalgebra whose elements are $ad$-semisimple and which
satisfies $\left[  \frak{s},\frak{t}\right]  =0$. 

\begin{proposition}
\textit{If }$\frak{g}=\frak{n}\oplus \frak{t}$ \textit{is rigid, then }$%
\frak{t}$ \textit{is a maximal torus.}
\end{proposition}

A proof can be found in [2].\newline Now, the interest of this result in connection with graph theory is the following 

\begin{theorem}
There do not exist decomposable, 2-step solvable rigid Lie algebras $\frak{r}=\frak{g}\oplus T$. 
\end{theorem}

\begin{proof}
If the algebra $\frak{r}$ is rigid, then the torus $T$ must be a
maximal torus of derivations for $\frak{g}$, and by proposition 4, the
algebra $\frak{r}$ is at least three step solvable, which contradicts the
assumption.
\end{proof}

\end{document}